\documentclass[12pt,twoside]{article}
\usepackage{amsmath}
\usepackage{mathrsfs}
\usepackage{amsfonts}
\usepackage{amssymb}

\usepackage{amsfonts,amsmath,amsthm}
\allowdisplaybreaks

\numberwithin{equation}{section} \textwidth 15.5 true cm
\begin{document}
  \title{{\bf Global existence and large time behavior for the system of compressible adiabatic flow through porous media in $\mathbb{R}^{3}$ }\footnote{Supported by National Natural Science Foundation
of China-NSAF (Grant No. 10976026).}}
  \author{ Guochun Wu \footnote{Corresponding author: Guochun Wu, Email
 address: guochunwu@126.com}\ \ , Zhong Tan and Jun Huang\\{\it\small School of Mathematical Sciences, Xiamen University,
 Fujian 361005, China}}
  \date{}
  \maketitle
  \begin{abstract}{\small
    The system of compressible adiabatic flow through porous media is considered in $\mathbb{R}^{3}$
    in the present paper. The global existence and uniqueness of classical solutions are obtained when the
    initial data is near its equilibrium. We also show that the pressure of the system converges to its equilibrium state at the same $L^2$-rate $(1+t)^{-\frac{3}{4}}$
    as the Navier-Stokes equations without heat conductivity, but the velocity of the system decays at the $L^2$-rate $(1+t)^{-\frac{5}{4}}$,
    which is faster than the $L^2$-rate $(1+t)^{-\frac{3}{4}}$ for the Navier-Stokes equations without heat conductivity [3]. }
  \end{abstract}

\section{Introduction}

$\ \ \ \ $The motion of compressible adiabatic flow can be modeled
by the compressible Euler equations with frictional damping terms,
say, the following balance laws:
$$
\left\{
\begin{array}{l}
\partial_t \rho+\nabla\cdot (\rho u)=0,\\
\partial_t (\rho u)+\nabla\cdot(\rho u\otimes u)+\nabla p=-a\rho u,\\
\partial_t (\rho\mathscr{E})+\nabla\cdot(\rho u\mathscr{E}+up)=-a\rho
u^2,
\end{array}
\right.\eqno (1.1)
$$
for $(t,x)\in[0,+\infty)\times\mathbb{R}^{3}$. Here $\rho$,
$u=(u_{1},u_{2},u_{3})^{t}$ and $p$ denote the density, the
velocity, and the pressure respectively. The total energy
$\mathscr{E}=\frac{|u|^2}{2}+e$, where $e$ is the internal energy.
The constant $a>0$ models friction. In this paper, we will consider
only polytropic fluids, so that the equations of state for the fluid
is given by
$$p=R\rho\theta,\ \ e=c_v\theta,$$
where $R > 0,\  c_v > 0$ are the universal gas constant and the
speci¡¥c heat at constant volume respectively.

In the case of isentropic flow where s=const., (1.1) takes the form
$$
\left\{
\begin{array}{l}
\partial_t \rho+\nabla\cdot (\rho u)=0,\\
\partial_t (\rho u)+\nabla\cdot(\rho u\otimes u)+\nabla p=-a\rho u.\\
\end{array}
\right.\eqno(1.2)
$$
For one dimension, the system (1.2) is hyperbolic with two
characteristic speeds $\lambda=u\pm\sqrt{p'(\rho)}$. As a vacuum
appears, it fails to be strict hyperbolic. Thus, the system involves
three mechanisms: nonlinear convection, lower-order dissipation of
damping and the resonance due to vacuum. The global existence of a
smooth solution with small data was first proved by Nishida [23,24],
and the behavior of the smooth solution was studied in many papers,
the reader can refer to [1,2,5-7,10,12-15,18,20,21,30,32] and
references therein. In multi-dimension space, Wang-Yang [29]
considered the time-asymptotic behavior of solutions in
multi-dimensions, the global existence and pointwise estimates of
the solutions were obtained. Sideris-Thomases-Wang [27] proved that
damping prevents the development of singularities in small amplitude
classical solutions in three-dimensional space, using an equivalent
reformulation of the Cauchy problem to obtain effective energy
estimates. Fang-Xu [4] studied the existence and asymptotic behavior
of $C^{1}$ solutions to the multi-dimensional space on the framework
of Besov space. Jang-Masmoudi [16] studied well-posedness of
compressible Euler equations in a physical vacuum. The optimal
estimates was obtained by Tan-Wu [28]. For initial boundary value
problem, refer for instance to [10,26] and references therein.

For the adiabatic flow where $s\neq$ const., much less is known. For
one dimension, the global existence of smooth solutions to the
Cauchy problem has been proved in [11] and [31] for small initial
data. The problem of large time behavior of these solutions is known
only for some particular initial data; see [8,22]. For initial
boundary value problem, refer for instance to [9,25] and references
therein.

From the physical point of view, the multi-dimensional case of model
(1.1) describes more realistic phenomena. Also the multi-dimensional
compressible Euler equations carry some unique features, such as the
effect of vorticity, which are totally absent in the one dimensional
case and make the problem more mathematical challenge, system (1.1)
and its time-asymptotic behavior are of great importance and are
much less understood than its one dimensional companion. To our
knowledge, there is not any work for the full system (1.1) in
$\mathbb R^3$.

The study of this paper is motivated by Duan-Ma [3], where the
authors considered the compressible Navier-Stokes equations without
heat conductivity in terms of the variables $p,u$ and $s$. It is
well-known that all thermodynamics variables $\rho,\theta,e,p$ as
well as the entropy $s$ can be denoted by functions of any two of
them. We take the two variables to be $p$ and $s$. Then the equation
of state for the gas is then given by
$$\rho=kp^{\frac{c_v}{c_v+R}}e^{-\frac{s}{c_v+R}},\eqno(1.3)$$
where $k>0$ is a constant. Under the aforementioned assumptions, the
system (1.1) in terms of the variables $p,u$ and $s$ reads
$$
\left\{
\begin{array}{l}
\partial_t p+\frac{R+c_v}{c_v}p\nabla\cdot u+u\cdot\nabla p=0,\\
\partial_t u+(u\cdot\nabla)u+\frac{\nabla p}{\rho}=-a u,\\
\partial_t s+(u\cdot\nabla)s=0,
\end{array}
\right.\eqno (1.4)
$$
where $\rho=\rho(p,s)$ is defined by (1.3). Notice that (1.4) is a
hyperbolic system, where the dissipation comes from damping. We
consider the initial value problem to (1.3) in the whole space
$\mathbb R^3$ with the initial data
$$(p(x),u(x),s(x))|_{t=0}=(p_0(x),u_0(x),s_0(x))\longrightarrow (p_{\infty},0,s_{\infty})\ \ as\ \ |x|\longrightarrow \infty,\eqno(1.5)$$
where $p_{\infty}>0$ and $s_{\infty}$ are given constants.

Before we state the main results, let us introduce some notations
for the use throughout this paper. $C$ denotes some positive
constant. The norms in the Sobolev Spaces $H^m(\mathbb R^3)$ and
$W^{m,q}(\mathbb R^3)$ are denoted respectively by $||\cdot||_m$ and
$||\cdot||_{m,q}$ for $m\geq 0$ and $q\geq 1$. In particular, for
$m=0$ we will simply use $||\cdot||$ and $||\cdot||_{L^q}$.
Moreover, we use $\langle\cdot ,\cdot\rangle$ to denote the inner
product in $L^2(\mathbb R^3)$. Finally,
$$\nabla=(\partial_{1},\partial_{2},\partial_{3}),\ \ \ \partial_i=\partial_{x_i},\ i=1,2,3,$$
and for any integer $l\geq 0$, $\nabla ^lf$ denotes all derivatives
of order $l$ of the function $f$. And for multi-indices $\alpha$ and
$\beta$ $$\alpha=(\alpha_1,\alpha_2,\alpha_3),\ \
\beta=(\beta_1,\beta_2,\beta_3),$$ we use
$$\partial_x^\alpha=\partial_{x_1}^{\alpha_1}\partial_{x_2}^{\alpha_2}\partial_{x_3}^{\alpha_3},\ \ |\alpha|=\sum_{i=1}^3\alpha_i,$$
and $C_\alpha^\beta=\frac{\alpha!}{\beta!(\alpha-\beta)!}$ where
$\beta\leq \alpha$.

For the global existence and large time behavior of classical
solutions, we have the following:\\
\\{\bf Theorem 1.1.}\ \ Let the initial data $(p_0,u_0,s_0)$ be such
that $||(p_0-p_{\infty},u_0,s_0-s_{\infty})||_{3}$ is sufficiently
small and $||(p_0-p_{\infty},u_0)||_{L^1}$ is bounded. Then the
initial value problem (1.4)-(1.5) admits a unique solution $(p,u,s)$
globally in time with $p>0$, satisfying
$$p-p_\infty,u,s-s_\infty\in C^0(0,\infty;H^3(\mathbb R^3))\cap C^1(0,\infty;H^2(\mathbb R^3)).$$
Moreover, there exists a constant $C_0$ such that for any $t\geq 0$,
$$||(p-p_\infty,u)(t)||^2_{3}+\int_0^t(||\nabla p(\tau)||_{2}^2+|| u(\tau)||_{3}^2)d\tau\leq C_0 ||(p_0-p_\infty,u_0)||_3^2,$$
$$||(s-s_\infty)(t)||_3\leq C_0||(p_0-p_\infty,u_0,s-s_\infty)||_3e^{C_0||(p_0-p_\infty,u_0)||_{L^1\cap H^3}}.$$
Finally, there is a constant $C_1$ such that for any $t\geq 0$, the
solution $(p,u,s)$ has the decay properties
$$
\begin{array}{l}
||(p-p_{\infty})(t)|| \leq C_1 (1+t)^{-\frac{3}{4}},\\
||\nabla (p-p_{\infty})(t)||_2+||u||_3\leq C_1(1+t)^{-\frac{5}{4}},\\
||\partial_t (p,u,s)(t)||\leq C_1(1+t)^{-\frac{5}{4}}.
\end{array}
$$

The rest of the paper is organized as follows. We will reformulate
the problem in Section 2. In section 3, we take Hodge decomposition
to analyze linearized system and establish the $L^2$ time-decay
rate. The proof of Theorem 1.1 is given in Section 4 and 5.

\section{Reformulated system}

$\ \ \ \ $ In this section, we first reformulate the problem as
follows. Set $$\kappa_1=\sqrt{\frac{c_v}{(R+c_v)\rho_\infty
p_\infty}},\ \kappa_2=\sqrt{\frac{(R+c_v)
p_\infty}{c_v\rho_\infty}},$$ where
$\rho_\infty=\rho(p_\infty,s_\infty)$. Taking change of variables by
$$(p,u,s)\longrightarrow (p+p_{\infty},\kappa_1u,s+s_\infty),$$
the initial value problem (1.4)-(1.5) is reformulated as
$$
\left\{
\begin{array}{l}
\partial_t p+\kappa_2\nabla\cdot u=F,\\
\partial_t u+\kappa_2\nabla p+au=G,\\
\partial_ts+\kappa_1(u\cdot\nabla)s=0,\\
(p,u,s)|_{t=0}:=(p_0,u_0,s_0)\rightarrow (0,0,0)\ as\ |x|\rightarrow
\infty,
\end{array}
\right.\eqno (2.1)
$$
where
$$
\begin{array}{rl}
F(p,u,s)&=-\frac{(R+c_v)\kappa_1}{c_v}p\nabla\cdot u-\kappa_1
u\cdot\nabla p,\\
G(p,u,s)&=-\kappa_1(u\cdot
\nabla)u-\frac{1}{\kappa_1}(\frac{1}{\rho}-\frac{1}{\rho_\infty})\nabla
p,.
\end{array}$$
Here and in the sequel, for the notational simplicity, we still
denote the reformulated variables by $(p,u,s)$.

Let us define the solution space and the solution norm of the
initial value problem (2.1) by
$$
\begin{array}{rl}
X(0,T)=&\{(p,u,s);p,u,s\in C^0(0,T;H^3(\mathbb R^3))\cap
C^1(0,T;H^2(\mathbb R^3)),\\&\nabla p\in L^2(0,T;H^2(\mathbb
R^3)),u\in L^2(0,T;H^3(\mathbb R^3))\},
\end{array}$$
and $$N(0,T)^2=\sup_{0\leq t\leq T}||(p,u,s)||_3^2+\int_0^T(||\nabla
p(t)||_2^2+||u(t)||_3^2)dt,$$ for any $0\leq T\leq \infty$. The
local existence and uniqueness of $H^3$ solution can be established
by following the methods in Kato [17] or Majda [19].\\
\\{\bf Proposition 2.1} (local existence){\bf .} Let $(p_0,u_0,s_0)\in H^3(\mathbb
R^3)$ be such that $$\inf_{x\in \mathbb R^3}\{p_0(x)+p_\infty\}>0.$$
Then there exists a positive constant $T_0$ depending on $N(0,0)$
such that the initial value problem (2.1) has a unique solution
$(p,u,s)\in X(0,T_0)$ which satisfies $N(0,T_0)\leq 2N(0,0)$ and
$$\inf_{x\in \mathbb R^3,0\leq t\leq T_0}\{p(t,x)+p_\infty\}>0.$$

To prove global existence of a smooth solution with small initial
data, we establish global a priori estimates of the solution.\\
\\{\bf Proposition 2.2} (A priori estimate){\bf .} Let $(p_0,u_0,s_0)\in H^3(\mathbb
R^3)$ and $(p_0,u_0)\in L^1(\mathbb R^3)$. Suppose that the initial
value problem (2.1) has a solution $(p,u,s)\in X(0,T)$, where T is a
positive constant. Then there exist a small constant $\epsilon>0$
and a constant $C_2$, which are independent of $T$, such that if
$$\sup_{0\leq t\leq T}||(p,u,s)(t)||_3\leq \epsilon,\eqno(2.2)$$
then for any $t\in [0,T]$, it holds that
$$||(p,u)(t)||_3^2+\int_0^t(||\nabla p(\tau)||_2^2+||u(\tau)||_3^2)d\tau\leq C_2||(p_0,u_0)||_3^2,\eqno(2.3)$$
$$||s(t)||_3\leq C_2||(p_0,u_0,s_0)||_3e^{C_2K_0},\eqno(2.4)$$
where $$K_0=||(p_0,u_0)||_{L^1\cap H^3}.\eqno(2.5)$$ Furthermore,
there is a constant $C_3$ such that for any $t\in [0,T]$, the
solution $(p,u,s)$ has the decay properties
$$
||p(t)|| \leq C_1 (1+t)^{-\frac{3}{4}},\eqno(2.6)$$
$$||\nabla p(t)||_2+||u||_3\leq C_1(1+t)^{-\frac{5}{4}},\eqno(2.7)$$
$$||\partial_t (p,u,s)(t)||\leq C_1(1+t)^{-\frac{5}{4}}.\eqno(2.8)$$

Theorem 1.1 follows from Proposition 2.1 and Proposition 2.2 by
standard continuity argument.

\section{Spectral analysis and linear $L^2$ estimates}

$\ \ \ \ $The linearized equations corresponding to system (2.1) is
$$
\left\{
\begin{array}{l}
\partial_t p+\kappa_2\nabla\cdot u=0,\\
\partial_t u+\kappa_2\nabla p+au=0,\\
\partial_ts=0,\\
(p,u,s)|_{t=0}:=(p_0,u_0,s_0)\rightarrow (0,0,0)\ as\ |x|\rightarrow
\infty.
\end{array}
\right.\eqno (3.1)
$$
We take Hodge decomposition to analyze the first two of linearized
system (3.1). For $r\in\mathbb R$, we denote by $\Lambda^{r}$ the
pseudo differential operator defined by
$\Lambda^{r}f=\mathscr{F}^{-1}(|\xi|^{r}\hat{f}(\xi))$. Let
$v=\Lambda^{-1}div\ u$ be the ``compressible part" of the velocity
and $\omega=\Lambda^{-1}curl\ u$ (with $(curl\
z)_i^j=\partial_jz^i-\partial_iz^j$) be the ``incompressible part",
then the first two of linearized system (3.1) writes
$$
\left\{
\begin{array}{l}
\partial_t p+\kappa_2\Lambda v=0,\\
\partial_t v-\kappa_2\Lambda p+av=0,\\
\partial_tw+aw=0.\\
\end{array}
\right.\eqno (3.2)
$$
Indeed, as the definition of $v$ and $\omega$, and relation
$$u=-\Lambda^{-1}\nabla v-\Lambda^{-1}div\  \omega\eqno(3.3)$$
involve pseudo-differential operators of degree zero, the estimates
in space $H^{l}(\mathbb{R}^{3})$ for the original function $u$ will
be the same as for $(v,\omega)$.

This section is devoted to the proof of the following results.\\
\\{\bf Proposition 3.1.} Let $U_0=(p_0,u_0)\in
H^{l}(\mathbb{R}^{3})\cap L^1(\mathbb{R}^{3})$, and $U=(p,u)$
satisfies the first two equation of system (3.1). Let
$v:=\Lambda^{-1}div\ u$ and $\omega:=\Lambda^{-1}curl\ u$. Then
there exists a constant $C$ such that for $0\leq |k|\leq l$,
$$
\|\partial_x^{k}\omega(t)\|\leq Ce^{-at}\|U_0\|_k,\eqno(3.4)
$$
$$
\|\partial_x^{k}p(t)\|\leq C(1+t)^{-\frac{3}{4}-\frac{|k|}{2}}
(\|U_0\|_{L^1}+ \|U_0\|_k),\eqno(3.5)
$$
$$
\begin{array}{ll}
&\|\partial_x^{k}v(t)\|+
\|\partial_x^{k}u(t)\|\\
\leq &C(1+t)^{-\frac{5}{4}-\frac{|k|}{2}} (\|U_0\|_{L^1}+
\|U_0\|_k).
\end{array}
\eqno(3.6)
$$
\\{\bf Proof:}\ \ The estimate for $\omega$ is obvious, so let us
focus on the first two equations of (3.2). In terms of the semigroup
theory for evolutionary equation, the solution $(p,v)$ of the first
two equation of system (3.2) can be expressed via the Cauchy problem
for $V=(p,v)^t$ as
$$V_t=BV,\ \ V(0)=V_0,\ \ t\geq 0.$$
which gives rise to $$V(t)=S(t)V_0=:e^{tB}V_0,\ \ t\geq 0.$$

What is left is to analyze the differential operator $B$ in terms of
its Fourier expression A and show the long time properties of the
semigroup $S(t)$. Taking the Fourier transform with respect to the
space variable yields
$$\frac{d}{dt}\hat{V}
=A(\xi)\hat{V}\ \ \hbox{with}\ \
A(\xi):=\left(\begin{array}{cc}0&-\kappa_2|\xi|\\\kappa_2|\xi|&-a\end{array}\right).
$$
The characteristic polynomial of $A(\xi)$ is
$\lambda^{2}+a\lambda+\kappa_2^{2}|\xi|^{2}$ and possesses two
distinct roots:
$$\lambda_{\pm}(\xi)=-\frac{a}{2}(1\pm \sqrt{1-\frac{4\kappa_2^{2}|\xi|^{2}}{a^{2}}}).\eqno(3.7)$$
The semigroup $e^{tA}$ is expressed as
$$e^{tA}=e^{\lambda_+t}P_++e^{\lambda_-t}P_-,$$
where the project operators $P_{\pm}$ can be computed as
$$P_+=\frac{A(\xi)-\lambda_-I}{\lambda_+-\lambda_-}\ \ \ and\ \ \ P_-=\frac{A(\xi)-\lambda_+I}{\lambda_--\lambda_+}.$$
By a direct computation, we can verify the exact expression the
Fourier transform $\hat G(\xi,t)$ of Green's function
$G(x,t)=e^{tB}$ as
$$
\begin{array}{rl}
\hat G(\xi,t)&=:e^{tA}=e^{\lambda_+t}P_++e^{\lambda_-t}P_-\\
&=\left(\begin{array}{cc}\frac{\lambda_+e^{\lambda_-t}-\lambda_-e^{\lambda_+t}}{\lambda_+-\lambda_-}&-\frac{\kappa_2|\xi|(e^{\lambda_+t}-e^{\lambda_-t})}{\lambda_+-\lambda_-}
\\\frac{\kappa_2|\xi|(e^{\lambda_+t}-e^{\lambda_-t})}{\lambda_+-\lambda_-}&\frac{\lambda_+e^{\lambda_-t}-\lambda_-e^{\lambda_+t}}{\lambda_+-\lambda_-}-\frac{a(e^{\lambda_+t}-e^{\lambda_-t})}{\lambda_+-\lambda_-}\end{array}\right)
.\end{array} \eqno(3.8)
$$

To derive the long-time decay rate of solutions whatever in $L^2$
framework or on point-wise estimate, we need to verify the
approximation of the eigenvalues (3.7), we are able to obtain that
it holds for $|\xi|\ll 1$ that
$$\frac{\lambda_+e^{\lambda_-t}-\lambda_-e^{\lambda_+t}}{\lambda_+-\lambda_-}\sim
\frac{-ae^{-\frac{\kappa_2^2|\xi|^2}{a^2}t}+\frac{\kappa_2^2|\xi|^2}{a^2}e^{-at}}{-a},\
\ |\xi|\ll 1,\eqno(3.9)$$
$$\frac{e^{\lambda_+t}-e^{\lambda_-t}}{\lambda_+-\lambda_-}\sim
\frac{e^{-at}-e^{-\frac{\kappa_2^2|\xi|^2}{a^2}t}}{-a},\ \ |\xi|\ll
1,\eqno(3.10)$$
$$\frac{\lambda_+e^{\lambda_-t}-\lambda_-e^{\lambda_+t}}{\lambda_+-\lambda_-}-\frac{a(e^{\lambda_+t}-e^{\lambda_-t})}{\lambda_+-\lambda_-}
\sim(1-\frac{\kappa_2^2|\xi|^2}{a^3})e^{-at},\ \ |\xi|\ll
1.\eqno(3.11)$$

To enclose the esimates, we also need to deal with the high
frequency $|\xi|\gg  1$. By a direct computation, we have the
following lemma:\\
\\{\bf Lemma 3.2.} For any given constant $\eta>0$, there exist two
positive constant $C$ and $R_0$ which are just dependent of $\eta$
such that $$|||\hat G(\xi,t)|||\leq Ce^{-R_0t}\ for \ \
|\xi|\geq\eta,$$ where
$|||\cdot|||$ denotes the norm of matrix for simplicity.\\

Let us now tackle the proof of (3.5) and (3.6). The proof relies on
the use of explicit expression for $\hat p$, $\hat v$. By (3.8), we
have
$$\hat p=\frac{\lambda_+e^{\lambda_-t}-\lambda_-e^{\lambda_+t}}{\lambda_+-\lambda_-}\hat p_0-\frac{\kappa_2|\xi|(e^{\lambda_+t}-e^{\lambda_-t})}{\lambda_+-\lambda_-}\hat v_0.$$
Combining (3.9), (3.10) and Lemma 3.2, we have the $L^2$-decay rate
for $p$ as
$$
\begin{array}{rl}
\|\hat{p}(t)\|_{L^{2}(\mathbb{R}^{3})}^{2}=
&\int_{|\xi|\leq\eta}|\hat{p}(\xi,t)|^{2}d\xi+
\int_{|\xi|\geq\eta}|\hat{p}(\xi,t)|^{2}d\xi\\
\leq&
C\int_{|\xi|\leq\eta}e^{-\frac{\kappa_2^{2}|\xi|^{2}}{a^2}t}(|\hat{p}_0|^{2}+|\hat{v}_0|^{2})d\xi\\
&+Ce^{-R_0t}\int_{|\xi|\geq\eta}(|\hat{p}_0|^{2}+|\hat{v}_0|^{2})d\xi\\
\leq& C\|(p_0,u_0)\|_{L^1(\mathbb{R}^{3})}^{2}
\int_{|\xi|\leq\eta}e^{-\frac{\kappa_2^{2}|\xi|^{2}}{a^2}t}d\xi\\
&+Ce^{-R_0t}\|(p_0,u_0)\|_{L^{2}(\mathbb{R}^{3})}^{2} \\
\leq &C(1+t)^{-\frac{3}{2}}\|(p_0,u_0)\|
_{L^{2}(\mathbb{R}^{3})\bigcap L^{1}(\mathbb{R}^{3})}^{2}.
\end{array}\eqno(3.12)
$$
Taking the same argument, we have the $L^2$-decay rate for $v$ as
$$
\begin{array}{rl}
\|\hat{v}(t)\|_{L^{2}(\mathbb{R}^{3})}^{2}=
&\int_{|\xi|\leq\eta}|\hat{v}(\xi,t)|^{2}d\xi+
\int_{|\xi|\geq\eta}|\hat{v}(\xi,t)|^{2}d\xi\\
\leq&
C\int_{|\xi|\leq\eta}e^{-\frac{\kappa_2^{2}|\xi|^{2}}{a^2}t}|\xi|^{2}(|\hat{p}_0|^{2}+|\hat{v}_0|^{2})d\xi\\
&+Ce^{-R_0t}\int_{|\xi|\geq\eta}(|\hat{p}_0|^{2}+|\hat{v}_0|^{2})d\xi\\
\leq &C(1+t)^{-\frac{5}{2}}\|(p_0,u_0)\|
_{L^{2}(\mathbb{R}^{3})\bigcap L^{1}(\mathbb{R}^{3})}^{2}.
\end{array}\eqno(3.13)
$$
The $L^{2}$-decay rate on the derivatives of $(n,u)$ as
$$
\begin{array}{rl}
\|\widehat{\partial_x^{k}p}(t)\|_{L^{2}(\mathbb{R}^{3})}^{2}=
&\int_{|\xi|\leq\eta}|\xi|^{2k}|\hat{p}(\xi,t)|^{2}d\xi+
\int_{|\xi|\geq\eta}|\xi|^{2k}|\hat{p}(\xi,t)|^{2}d\xi\\
\leq&
C\int_{|\xi|\leq\eta}|\xi|^{2k}e^{-\frac{\kappa_2^{2}|\xi|^{2}}{a^2}t}(|\hat{p}_0|^{2}+|\hat{v}_0|^{2})d\xi\\
&+Ce^{-R_0t}\int_{|\xi|\geq\eta}|\xi|^{2k}(|\hat{p}_0|^{2}+|\hat{v}_0|^{2})d\xi\\
\leq &C(1+t)^{-\frac{3}{2}-k}\|(p_0,u_0)\|
_{H^{k}(\mathbb{R}^{3})\bigcap L^{1}(\mathbb{R}^{3})}^{2},
\end{array}\eqno(3.14)
$$
and
$$
\begin{array}{rl}
\|\widehat{\partial_x^{k}v}(t)\|_{L^{2}(\mathbb{R}^{3})}^{2}=
&\int_{|\xi|\leq\eta}|\xi|^{2k}|\hat{v}(\xi,t)|^{2}d\xi+
\int_{|\xi|\geq\eta}|\xi|^{2k}|\hat{v}(\xi,t)|^{2}d\xi\\
\leq&
C\int_{|\xi|\leq\eta}|\xi|^{2k}e^{-\frac{\kappa_2^{2}|\xi|^{2}}{a^2}t}|\xi|^{2}(|\hat{p}_0|^{2}+|\hat{v}_0|^{2})d\xi\\
&+Ce^{-R_0t}\int_{|\xi|\geq\eta}|\xi|^{2k}(|\hat{p}_0|^{2}+|\hat{v}_0|^{2})d\xi\\
\leq &C(1+t)^{-\frac{5}{2}-k}\|(p_0,u_0)\|
_{L^{2}(\mathbb{R}^{3})\bigcap L^{1}(\mathbb{R}^{3})}^{2},
\end{array}\eqno(3.15)
$$
for $k\geq 1$. By the relation between $u$ and $(v,\omega)$, we can
easily get the estimate for $u$:
$$\|\partial_x^{k}u(t)\|\leq C(1+t)^{-\frac{5}{2}-k}
\|U_0\|_{H^{k}(\mathbb{R}^{3})\bigcap L^{1}(\mathbb{R}^{3})}.$$ The
proof of proposition 3.1 is
completed.\\

We also need the following Sobolev's inequalities.\\
\\{\bf Lemma 3.3.} Let $f\in H^2(\mathbb R^3)$. Then it holds:\\
$(i)\ \ \ ||f||_{L^\infty}\leq C||\nabla f||^{\frac{1}{2}}||\nabla
f||_1^{\frac{1}{2}} \leq C||\nabla f||_1$;\\
$(ii)\ \ ||f||_{L^6}\leq C||\nabla f||$;\\
$(iii)\ ||f||_{L^q}\leq C||f||_1,\ 2\leq q\leq 6$.

\section{A priori estimates}

$\ \ \ \ $We suppose that the inequality (2.2) holds throughout this
section and next section. The initial value problem (2.1) has a
solution $(p,u,s)$ in the space $X(0,T)$ with some $T\in
(0,+\infty]$. We also omit the variable $t$ of all functions in the
proof of different lemmas in this section for simplicity.

In what follows, a series of lemmas on the energy estimates are
given. Firstly, the energy estimate of lower order for $(p,u)$ is
obtained in the following lemma. \\
\\{\bf Lemma 4.1.} There exists a constant $D_1>0$ suitably large
which is independent of $\epsilon$ such that
$$\frac{d}{dt}(D_1||(p,u)(t)||^2+\langle \nabla p,u\rangle (t))+C(||\nabla p(t)||^2+||u(t)||^2)\leq C||\nabla u(t)||^2,\eqno(4.1)$$
for any $0\leq t\leq T$.\\
\\{\bf Proof:}\ \ Multiplying $(2.1)_1$-$(2.1)_2$ by $p,\ u$
respectively and then integrating them over $\mathbb R^3$, we have
$$\frac{1}{2}\frac{d}{dt}||(p,u)(t)||^2+a||u(t)||^2=\langle p,F \rangle +\langle u,G \rangle.\eqno(4.2)$$
The two terms on the right hand side of the above equation can be
estimated as follows.

Firstly, for the first term, it holds that
$$\begin{array}{rl}\langle p,F
\rangle&=-\frac{(R+c_v)\kappa_1}{c_v}\langle p,p\nabla\cdot
u\rangle-\kappa_1 \langle p,u\cdot\nabla
p\rangle\\&=(-\frac{(R+c_v)\kappa_1}{c_v}+\frac{\kappa_1}{2})\langle
p^2,\nabla\cdot u\rangle.\end{array}\eqno (4.3)$$ It follows from
Lemma 3.3, H\"older's inequality and (2.2) that
$$\begin{array}{rl}|\langle
p^2,\nabla\cdot u\rangle|&\leq ||p||_{L^3}||p||_{L^6}||\nabla
u||\leq C||p||_1||\nabla p||||\nabla u||\\&\leq C\epsilon (||\nabla
p||^2+||\nabla u||^2).\end{array}\eqno(4.4)$$  Putting (4.4) into
(4.3,) we arrive at
$$|\langle p,F \rangle|\leq C\epsilon(||\nabla
p||^2+||\nabla u||^2+||u||^2).\eqno(4.5)$$ For the second term, we
have$$|\langle u,G \rangle|\leq C(|\langle u,(u\cdot\nabla)u
\rangle|+|\langle u,(\frac{1}{\rho}-\frac{1}{\rho_\infty})\nabla p
\rangle|).\eqno(4.6)$$ Similar to the proof of (4.5), it follows
from Lemma 3.3, H$\ddot{o}$lder inequality and (2.2) that
$$|\langle u,(u\cdot\nabla)u
\rangle|\leq ||u||_{L^3}||u||_{L^6}||\nabla u||\leq C||u||_1||\nabla
u||^2\leq C\epsilon ||\nabla u||^2,\eqno(4.7)$$
$$\begin{array}{rl}|\langle u,(\frac{1}{\rho}-\frac{1}{\rho_\infty})\nabla p
\rangle|&\leq
||u||_{L^6}||\frac{1}{\rho}-\frac{1}{\rho_\infty}||_{L^3}||\nabla
p||\\&\leq C||\nabla u||||(p,s)||_1||\nabla p||\\&\leq
C\epsilon(||\nabla p||^2+||\nabla u||^2),\end{array}\eqno(4.8)$$
where by (1.3) and (2.2), we have used the fact
$$\rho\sim\rho_\infty+\mathcal {O}(1)(p+s)$$
and$$\frac{1}{\rho}-\frac{1}{\rho_\infty}\sim\mathcal {O}(1)(p+s).$$
Substituting (4.7) and (4.8) into (4.6), we obtain that the second
term is bounded by
$$|\langle u,G \rangle|\leq C\epsilon(||\nabla p||^2+||\nabla u||^2).\eqno(4.9)$$
Hence combining (4.2), (4.5) and (4.9)
yields$$\frac{d}{dt}||(p,u)||^2+C||u||^2\leq C\epsilon(||\nabla
p||^2+||\nabla u||^2),\eqno(4.10)$$ since $\epsilon> 0$ is
sufficiently small.

Next we shall estimate $||\nabla p||^2$. From $(2.1)_2$, we have
$$\kappa_2||\nabla p||^2=\langle -u_t,\nabla p \rangle-a\langle u,\nabla p \rangle+\langle G,\nabla p \rangle.\eqno(4.11)$$
By $(2.1)_1$, the first term on the right hand side can be written
as $$\begin{array}{rl}\langle -u_t,\nabla p
\rangle&=-\frac{d}{dt}\langle \nabla p,u \rangle+\langle \nabla
p_t,u \rangle\\&=-\frac{d}{dt}\langle \nabla p,u \rangle-\langle
 p_t,\nabla\cdot u \rangle\\&=-\frac{d}{dt}\langle \nabla p,u \rangle+\langle
 \kappa_2\nabla\cdot u-F,\nabla\cdot u \rangle.\end{array}\eqno(4.12)$$
It follows from the definition of $F$ that
$$\begin{array}{rl}|\langle -F,\nabla\cdot u\rangle|&\leq C(|\langle p\nabla\cdot u,\nabla\cdot u\rangle|
+|\langle u\cdot\nabla p,\nabla\cdot u\rangle|)\\&\leq C(||\nabla
p||_1||\nabla\cdot u||^2+||\nabla p||_1||\nabla
u||^2+||u||_1||\nabla u||^2)\\&\leq C\epsilon ||\nabla
u||^2.\end{array}\eqno(4.13)$$ Taking the same argument to estimate
$\langle G,\nabla p\rangle$, we have
$$|\langle G,\nabla p\rangle|\leq C\epsilon(||\nabla p||^2+||\nabla u||^2).\eqno(4.14)$$
Using Cauchy-Schwarz inequality, we easily get
$$|a\langle u,\nabla p\rangle|\leq C||u||^2+\frac{\kappa_2}{4}||\nabla p||^2.\eqno(4.15)$$
Since $\epsilon>0$ is small enough, putting (4.12), (4.13) (4.14)
and (4.15) into (4.11) gives
$$\frac{d}{dt}\langle\nabla p,u\rangle +\frac{\kappa_2}{2}||\nabla p||^2\leq C||u||+C\epsilon||\nabla u||.\eqno(4.16)$$
Multiplying (4.10) by $D_1$ suitably large and adding it to (4.16),
we finally deduce the lemma since $\epsilon>0$ is sufficiently
small. This completes the proof of Lemma 4.1.\\

Our next goal is to deal the higher order estimate of $(p,u)$.\\
\\{\bf Lemma 4.2.} For  any $0\leq t\leq T$, there exists a constant $D_2>0$ sufficiently
large which is independent of $\epsilon$, $t$ and $T$ such that
$$\begin{array}{ll}\frac{d}{dt}\{D_2H_1(p(t),u(t))+\sum\limits_{1\leq|\alpha|\leq 2}\langle\partial_x^\alpha\nabla p
,\partial_x^{\alpha}u\rangle(t)\}\\+C(||\nabla^2p(t)||_1^2+||\nabla
u(t)||_2^2)\leq C\epsilon(||\nabla
p(t)||^2+||u(t)||^2),\end{array}\eqno(4.17)$$ where $H_1(p,u)$ is equivalent to $||\nabla (p,u)(t)||_2^2$, if $\epsilon$ is small enough.\\
\\{\bf Proof:}\ \ For each multi-index $\alpha$ with $1\leq |\alpha|\leq
3$, by applying $\partial_x^\alpha$ to $(2.1)_1-(2.1)_2$,
multiplying them by $\partial_x^\alpha p,\ \partial_x^\alpha u$
respectively and then integrating them over $\mathbb R^3$, we have
$$\begin{array}{rl}&\frac{1}{2}\frac{d}{dt}||\partial_x^\alpha(p,u)||^2+a||\partial_x^\alpha u||
\\=&\langle\partial_x^\alpha p,\partial_x^\alpha F\rangle+\langle\partial_x^\alpha u,\partial_x^\alpha G\rangle\\=&I_1+I_2,\end{array}\eqno(4.18)$$
where $I_i,\ i=1,2$ are the corresponding terms in the above
equation which will be estimated as follows. Here and in the sequel
proof, the repeated index denotes summation over the index.

Firstly, for $I_1$, it holds that,
$$\begin{array}{rl}I_1=&-\frac{(R+c_v)\kappa_1}{c_v}\langle \partial^\alpha_x p,\partial^\alpha_x(p\nabla\cdot u)\rangle
-\kappa_1\langle\partial^\alpha_x p,\partial^\alpha_x(u\cdot \nabla
p)\rangle\\=&-\frac{(R+c_v)\kappa_1}{c_v}\langle \partial^\alpha_x
p,p\partial^\alpha_x(\nabla\cdot u)\rangle
-\kappa_1\langle\partial^\alpha_x p,(u\cdot \nabla
\partial^\alpha_x p)\rangle\\&-\kappa_1\sum\limits_{|\beta|\leq
|\alpha|-1}\langle\partial^\alpha_x
p,\partial^{\alpha-\beta}_xu\cdot \nabla
\partial^{\beta}_xp\rangle\\&-\sum\limits_{|\beta|\leq
|\alpha|-1}\frac{(R+c_v)\kappa_1}{c_v}C_\alpha^\beta\langle
\partial^\alpha_x p,\partial^{\alpha-\beta}_xp\partial^\beta_x(\nabla\cdot
u)\rangle\\=&\sum\limits_{i=1}^4I_{1,i}.
\end{array}\eqno(4.19)$$
By symmetry, Lemma 3.3 and some tedious but straightforward
calculations, we can obtain$$\sum\limits_{i=2}^4I_i\leq
C||(p,u,s)||_3(||\nabla p||_2^2+||u||_3^2)\leq C\epsilon(||\nabla
p||_2^2+||u||_3^2).\eqno(4.20)$$ From $(2.1)_1$, we have
$$div\ u=-\frac{c_v}{\kappa_2c_v+(R+c_v)\kappa_1p}(p_t+\kappa_1u\cdot\nabla p+\frac{a\kappa_1^2\rho u^2}{c_v}),$$
$$p_t=-\frac{\kappa_2c_v+(R+c_v)\kappa_1p}{c_v}div\ u-\kappa_1u\cdot\nabla p-\frac{a\kappa_1^2\rho u^2}{c_v},$$
which imply that
$$\begin{array}{rl}I_{1,1}=&(R+c_v)\langle \partial^\alpha_x
p,p\partial^\alpha_x(\frac{1}{\kappa_2c_v+(R+c_v)\kappa_1p}(p_t+\kappa_1u\cdot\nabla
p+\frac{a\kappa_1^2\rho u^2}{c_v}))\rangle\\
\leq&\frac{(R+c_v)}{2}\langle ((\partial^\alpha_x
p)^2)_t,\frac{p}{\kappa_2c_v+(R+c_v)\kappa_1p}\rangle+C\epsilon(||\nabla
p||_2^2+||u||_3^2)\\\leq&\frac{(R+c_v)}{2}\frac{d}{dt}\langle
(\partial^\alpha_x
p)^2,\frac{p}{\kappa_2c_v+(R+c_v)\kappa_1p}\rangle+C\epsilon(||\nabla
p||_2^2+||u||_3^2).
\end{array}\eqno(4.21)$$
Combining (4.19) and (4.21), we have
$$I\leq\frac{(R+c_v)}{2}\frac{d}{dt}\langle (\partial^\alpha_x
p)^2,\frac{p}{\kappa_2c_v+(R+c_v)\kappa_1p}\rangle+C\epsilon(||\nabla
p||_2^2+||u||_3^2).$$ Taking the same argument to deal $I_2$, then
there exists a function $H_1(p,u)$ which is equivalent to $||(\nabla
p,\nabla u)||_2^2$ and satisfies
$$\frac{d}{dt}H_1(p(t),u(t))+||\nabla u||_2^2\leq C\epsilon(||\nabla p||_2^2)+||u||^2).\eqno(4.22)$$
By a direct computation, we have the estimate on
$||\nabla\partial^\alpha_x p||$ for $1\leq|\alpha|\leq2$ as
following
$$\frac{\kappa_1}{2}\sum\limits_{1\leq\alpha\leq2}||\nabla\partial^\alpha_x p||^2+\frac{d}{dt}\sum\limits_{1\leq|\alpha|\leq2}
\langle\partial^\alpha_xu,\nabla\partial^\alpha_xp\rangle\leq
C\epsilon(||\nabla p||_2^2)+||u||^2_3).\eqno(4.23)
$$
Since $\epsilon$ is sufficiently small, multiplying (4.22) by $D_2$
suitably large and adding it to (4.23) give (4.17). Thus we
completes the proof of the lemma.\\

Finally, by symmetry and some tedious but straightforward
calculation, we have the energy estimates on the entropy.\\
\\{\bf Lemma 4.3.} It holds that
$$\frac{d}{dt}||s(t)||\leq C||u(t)||_3||s(t)||^2_3,\eqno(4.24)$$
for any $0\leq t\leq T$.

\section{The proof of global well-posedness}

$\ \ \ \ $In this section, we are devoted to prove Proposition 2.2.
We first consider a priori decay-in-time estimates on $(\nabla
p,u)$. This will be based on Proposition 3.1 about the decay
estimates on the linearized system. The decay-in-time estimate on
$p$ can be derived by decay-in-time estimates on $(\nabla p,u)$.
Precisely, we have the
following lemma.\\
\\{\bf Lemma 5.1.} Let $(p,u,s)$ be the solution of (2.1), then
$p$ satisfies the following inequality
$$||\nabla p||\leq CK_0(1+t)^{-\frac{5}{4}}+C\epsilon\int_0^t(1+t-\tau)^{-\frac{5}{4}}(||\nabla p(\tau)||_2+||u(\tau)||_3)d\tau,\eqno(5.1)$$
for any $0\leq t\leq T$, where $K_0=||(p_0,u_0)||_{L^1\cap H^3}$ as
in (2.5).\\
\\{\bf Proof:} From the Duhamel's principle and Proposition 3.1, we
have$$||\nabla p||\leq
CK_0(1+t)^{-\frac{5}{4}}+C\epsilon\int_0^t(1+t-\tau)^{-\frac{5}{4}}(||(F,G)(\tau)||_{L^1\cap
H^1})d\tau.\eqno(5.2)$$ The nonlinear term $(F,G)$ can be estimated
as following$$\begin{array}{rl}||(F,G)||_{L^1}&\leq
C||(p,u,s)(t)||_1(||\nabla p(t)||_1+||u(t)||_2)\\&\leq
C\epsilon(||\nabla p(t)||_1+||u(t)||_2),\end{array}\eqno(5.3)$$
$$\begin{array}{rl}||(F,G)||_{1}&\leq
C||(p,u,s)(t)||_{W^{1,\infty}}(||\nabla p(t)||_2+||u(t)||_3)\\&\leq
C\epsilon(||\nabla p(t)||_2+||u(t)||_3),\end{array}\eqno(5.4)$$
Putting (5.3) and (5.4) into (5.2), we obtain (5.1). The proof of
Lemma 5.1 is completed.\\

Now we are in a position to prove Proposition 2.2.\\
\\{\bf Proof of Proposition 2.2.}  We do it by three steps.\\

{\bf Step 1:} Since $\epsilon>0$ is sufficiently small, from Lemma
4.1 and Lemma 4.2, we have a function $H_2(p,u)$ which is equivalent
to $||(p,u)||^2_3$ and satisfies
$$\frac{d}{dt}H_2(p(t),u(t))+C(||\nabla p(t)||_2^2+|| u||_3^2)\leq 0,\eqno(5.5)$$
for any $0\leq t\leq T$, which implies (2.3).\\

{\bf Step 2:} Multiplying $(2.1)_2$ by $u$, integrating over
$\mathbb R^3$ and using Cauchy-Schwarz inequality, we have
$$\frac{d}{dt}||u||^2+C||u||^2\leq C||\nabla p||^2+C\epsilon||\nabla u||^2.\eqno(5.6)$$
Now we define the temporal energy functional
$$H_3(t)=||u(t)||^2+D_2H_1(p(t),u(t))+\sum\limits_{1\leq|\alpha|\leq 2}\langle\partial_x^\alpha\nabla p
,\partial_x^{\alpha}u\rangle(t)$$ for any $0\leq t\leq T$, where it
is noticed that $H_3(t)$ is equivalent to $||\nabla
p||_2^2+||u||_3^2$ since $D_2$ can be large enough.

Combining Lemma 4.2 and (5.6), we obtain
$$\frac{d}{dt}H_3(t)+C(||\nabla^2 p(t)||_1^2+||u(t)||_3^2)\leq C||\nabla p||^2.$$
Adding $||\nabla p||^2$ to both sides of the above inequality gives
$$\frac{d}{dt}H_3(t)+D_3H_3(t)\leq C||\nabla p||^2,\eqno(5.7)$$
where $D_3$ is a positive constant independent of $\epsilon$.
Set$$M(t)=\sup_{0\leq\tau \leq
t}(1+\tau)^{\frac{5}{2}}H_3(\tau),\eqno(5.8)$$ and notice
that$$||\nabla p(\tau)||_2+||u||_3\leq C\sqrt{H_3(\tau)}\leq
C(1+\tau)^{-\frac{5}{4}}\sqrt{M(t)},\ 0\leq\tau\leq t\leq
T.\eqno(5.9)$$ Then it follows from Lemma 5.1 that
$$\begin{array}{rl}||\nabla p(t)||&\leq CK_0(1+t)^{-\frac{5}{4}}+C\epsilon\int_0^t(1+t-\tau)^{-\frac{5}{4}}(1+\tau)^{-\frac{5}{4}}d\tau\sqrt{M(t)}
\\&\leq C(1+t)^{-\frac{5}{4}}(K_0+\epsilon\sqrt{M(t)}).
\end{array}\eqno(5.10)$$
Hence, by the Gronwall's inequality, (5.7) and (5.10) lead to
$$\begin{array}{rl}H_3(t)&\leq e^{-D_3t}H_3(0)+C\int_0^te^{-D_3(t-\tau)}||\nabla p(\tau)||^2d\tau\\
&\leq
e^{-D_3t}H_3(0)+C\int_0^te^{-D_3(t-\tau)}(1+\tau)^{-\frac{5}{2}}d\tau(K_0^2+\epsilon^2M(t))\\
&\leq
C(1+t)^{-\frac{5}{2}}(K_0^2+\epsilon^2M(t)).\end{array}\eqno(5.11)$$
Since $M(t)$ is non-decreasing, we have from (5.8) and (5.11) that
$$M(t)\leq C(K_0^2+\epsilon^2M(t)),\eqno(5.12)$$
for any $0\leq t\leq T$, which implies that $$M(t)\leq
CK_0^2,\eqno(5.13)$$ since $\epsilon>0$ is small enough. Thus we
obtain (2.7) from (5.9) and (5.13).

Next, using Proposition 3.1, (5.3) and (5.4), it follows from the
Duhamel's principle that
$$\begin{array}{rl}||p(t)||&\leq CK_0(1+t)^{-\frac{3}{4}}+C\int_0^t(1+t-\tau)^{-\frac{3}{4}}||(F,G)(\tau)||_{L^1\cap
L^2}d\tau\\&\leq
CK_0(1+t)^{-\frac{3}{4}}+C\int_0^t(1+t-\tau)^{-\frac{3}{4}}(||\nabla
p(\tau)||_{2}+||u||_3)d\tau\\& \leq
CK_0(1+t)^{-\frac{3}{4}}+CK_0\int_0^t(1+t-\tau)^{-\frac{3}{4}}(1+\tau)^{-\frac{5}{4}}d\tau\\
&\leq K_0(1+t)^{-\frac{3}{4}},\end{array}$$ for any $0\leq t\leq T$.
Thus (2.6) is proved.\\

{\bf Step 3:} By Lemma 4.3 and Gronwall's inequality, we arrive at
$$\begin{array}{rl}s(t)&\leq s(0)\exp\{C\int_0^t||u(\tau)||_3d\tau\}\\
&\leq s(0)\exp \{CK_0\int_0^t(1+\tau)^{-\frac{5}{4}}d\tau\}\\&\leq
s(0)\exp\{CK_0\}
\end{array}$$
i.e., for any $0\leq t\leq T$,$$||(p,u,s)(t)||_3^2\leq
C||(p_0,u_0,s_0)||_3^2\exp\{CK_0\}.$$ Hence (2.4) holds. For (2.8),
using the above estimates and (2.1) we have
$$\begin{array}{rl}||\partial_t(p,u,s)(t)||&\leq C(||u||_1+||\nabla
p||)\\&\leq CK_0(1+t)^{-\frac{5}{4}}.
\end{array}$$
for any $0\leq t\leq T$. Thus, (2.8) is proved and this completes
the proof of Proposition 2.2.


\begin{thebibliography}{aa}
\footnotesize
 \bibitem{CW1}D. Coutand, S. Shkoller, Well-posedness in smooth function spaces for the moving-boundary 1-D compressible Euler equations in physical vacuum, Comm. Pure Appl. Math., 64  (2011), 328--366.






 \footnotesize \vspace{-0.3cm}
    \bibitem{DF}C. M. Dafermos, Can dissipation prevent the breaking of waves? In:
Transactions of the Twenty-Sixth Conference of Army Mathematicians,
187-198, ARO Rep. 81, 1, U. S. Army Res. Office, Research Triangle
Park, N.C., 1981.

\footnotesize \vspace{-0.3cm}
    \bibitem{RL}R. Duan, H. Ma, Global Existence and Convergence Rates for the
3-D Compressible Navier-Stokes Equations without Heat Conductivity,
Indiana Univ. Math. J. 57 (2008), 2299-2319.

\footnotesize \vspace{-0.3cm}
    \bibitem{RL}D. Y. Fang, J. Xu, Existence and asymptotic behavior
    of $C^1$ solutions to the multi-dimensional compressible Euler
    equations with damping, Nonlinear Analysis 70 (2009) 244-261.

\footnotesize \vspace{-0.3cm}
    \bibitem{DF}L. Hsiao, Quasilinear Hyperbolic Systems and Dissipative
Mechanisms, World Scientific, Singapore, 1998.


\footnotesize \vspace{-0.3cm}
    \bibitem{DF}L. Hsiao, T. P. Liu, Convergence to nonlinear diffusion waves for
solutions of a system of hyperbolic conservation laws with damping,
Comm Math Phys 143 (1992), 599-605.

\footnotesize \vspace{-0.3cm}
    \bibitem{DF}L. Hsiao, T. P. Liu, Nonlinear diffusive phenomena of nonlinear
hyperbolic systems, Chinese Ann. Math. Ser. B 14 (1993), 1-16.

\footnotesize \vspace{-0.3cm}
    \bibitem{DF}L. Hsiao, T. Luo, Nonlinear diffusive phenomena of
    solutions for the system of compressible adiabatic flow through
    porous media, J. Differential Equations 125 (1996), 329-365.



\footnotesize \vspace{-0.3cm}
    \bibitem{DF}L. Hsiao, R. H. Pan, Initial-boundary value problem for the system of compressible
adiabatic flowt hrough porous media. J. Differential Equations 159
(1999) 280-305.

\footnotesize \vspace{-0.3cm}
    \bibitem{DF}L. Hsiao, R. H. Pan, The damped p-system with boundary
    effects,
Contemporary Mathematics, 255(2000), 109-123.




 \footnotesize \vspace{-0.3cm}
    \bibitem{DF}L. Hsiao, D. Serre, Global existence of solutions
for the system of compressible adiabatic flow through porous media,
SIAM J. Math Anal. 27 (1996), 70-77.



\footnotesize \vspace{-0.3cm}
    \bibitem{DF}F. M. Huang, P. Marcati, R. H. Pan, Convergence to Barenblatt
Solution for the Compressible Euler Equations with Damping and
Vacuum, Arch. Ration. Mech. Anal. 176 (2005), 1-24.

\footnotesize \vspace{-0.3cm}
    \bibitem{DF}F. M. Huang, R. H. Pan, Asymptotic behavior of the solutions to
the damped compressible Euler equations with vacuum, J. Differential
Equations, 220 (2006), 207-233.

\footnotesize \vspace{-0.3cm}
    \bibitem{DF}F. M. Huang, R. H. Pan, Convergence rate for compressible Euler
equations with damping and vacuum, Arch. Ration. Mech. Anal. 166
(2003) 359-376.

\footnotesize \vspace{-0.3cm}
    \bibitem{DF}J. Jang, N. Masmoudi, Well-posedness for compressible Euler
equations with physical vacuum singularity,  Comm. Pure Appl. Math.
62 (2009), 1327-1385.

\footnotesize \vspace{-0.3cm}
    \bibitem{DF}J. Jang, N. Masmoudi, Well-posedness of compressible Euler equations in a physical vacuum
arXiv: 1005. 4441.






\footnotesize \vspace{-0.3cm}
    \bibitem{DF}T. Kato, The Cauchy problem for quasi-linear symmetric hyperbolic
systems, Arch. Rational Mech. Anal. 58 (1975) 181-205.

\footnotesize \vspace{-0.3cm}
    \bibitem{DF}T. P. Liu, Compressible flow with damping and vacuum, Japan J. Appl.
Math, 13 (1996), 25-32.






\footnotesize \vspace{-0.3cm}
    \bibitem{DF}A. Majda, Compressible Fluid Flow and Conservation laws in Several Space Variables,
    Springer-
Verlag, Berlin/New York, 1984.

\footnotesize \vspace{-0.3cm}
    \bibitem{DF}P. Marcati, A. Milani, The one-dimensional Darcy's law as the limit of a compressible Euler flow, J. Differential Equations 84 (1990), no. 1,
    129-147.

\footnotesize \vspace{-0.3cm}
    \bibitem{DF}P. Marcati, B. Rubino, Hyperbolic to Parabolic Relaxation Theory for Quasilinear First Order Systems, J. Differential Equations 162 (2000), no.2,
    359--399.

\footnotesize \vspace{-0.3cm}
    \bibitem{DF}P. Marcati, R. H. Pan, On the diffusive profiles for the system of compressible adiabatic flow through porous
    media, SIAM J. Math. Anal. 33 (2001), 790-826.



\footnotesize \vspace{-0.3cm}
    \bibitem{DF}T. Nishida, Global solutions for an initial-boundary
    value problem of a quasilinear hyperbolic systems, Proc. Japan
    Acad. 44 (1968) 642-646.

\footnotesize \vspace{-0.3cm}
    \bibitem{DF}T. Nishida, Nonlinear hyperbolic equations and
    relates topics in fluid dynamics, Publ. Math. D'Orsay (1978)
    46-53.

\footnotesize \vspace{-0.3cm}
    \bibitem{DF}R. H. Pan, Boundary effects and large time behavior
    for the system of compressible adiabatic flow through porous
    media, Michigan Math. J. 49 (2001), 519-539.



\footnotesize \vspace{-0.3cm}
    \bibitem{DF}R. H. Pan, K. Zhao, The 3D compressible Euler
    equations with damping in a bounded domain. J. Differential
    Equations 246 (2009) 581-596.


\footnotesize \vspace{-0.3cm}
    \bibitem{DF}T. C. Sideris, B. Thomases, D. H. Wang, Long time behavior of
solutions to the 3D compressible Euler equations with damping, Comm.
Partial Differential Equations 28 (2003) 795-816.

\footnotesize \vspace{-0.3cm}
    \bibitem{DF}Z. Tan, G. C. Wu, Large time behavior of solutions for compressible Euler equations with damping in
    $\mathbb{R}^{3}$, J. Differential Equations (in press).



\footnotesize \vspace{-0.3cm}
    \bibitem{DF}W. Wang, T. Yang, The pointwise estimates of solutions for Euler
equations with damping in multi-dimensions. J. Differential
Equations, 173 (2001), 410-450.


 \footnotesize \vspace{-0.3cm}
    \bibitem{DF}H. J. Zhao, Convergence to strong nonlinear diffusion waves for
solutions of p-system with damping, J. Differential Equations, 174 (
2001), 200-236.

 \footnotesize \vspace{-0.3cm}
    \bibitem{DF}Y. Zheng, Global smooth solutions to the adiabatic gas dynamics
system with dissipation terms, Chinese Ann. Math. 17A (1996),
155-162.


\footnotesize \vspace{-0.3cm}
    \bibitem{DF}C. J. Zhu, Convergence rates to nonlinear diffusion waves for weak
entropy solutions to p-system with damping, Science in China, Ser.
A, 46 (2003), 562-575.


  \end{thebibliography}
\end{document}